\documentclass[oneside,english,reqno]{amsart}
\usepackage[T1]{fontenc}
\usepackage[latin1]{inputenc}
\usepackage{graphicx}

\makeatletter
 \theoremstyle{plain}    
 \newtheorem{thm}{Theorem}[section]
 \numberwithin{equation}{section} 
 \numberwithin{figure}{section} 
 \theoremstyle{plain}
 \theoremstyle{remark}    
 \newtheorem*{acknowledgement*}{Acknowledgement} 
 \theoremstyle{definition}
 \newtheorem{defn}[thm]{Definition}
 \theoremstyle{plain}    
 \newtheorem{lem}[thm]{Lemma} 
 \theoremstyle{plain}    
 \newtheorem{prop}[thm]{Proposition} 

\DeclareMathOperator{\sd}{sd}
\DeclareMathOperator{\cone}{cone}
\DeclareMathOperator{\Bier}{Bier}
\DeclareMathOperator{\ub}{ub}
\DeclareMathOperator{\lb}{lb}
\usepackage{babel}
\makeatother
\begin{document}

\title{Nested set complexes for posets and the Bier construction}

\author{Juliane Lehmann}

\begin{abstract}
We generalize the concept of combinatorial nested set complexes to
posets and exhibit the topological relationship between the arising
nested set complexes and the order complex of the underlying poset.
In particular, a sufficient condition is given so that this relationship
is actually a subdivision.

We use the results to generalize the proof method of \v{C}uki\'{c}
and Delucchi, so far restricted to semilattices, for a result of Björner,
Paffenholz, Sjöstrand and Ziegler on the Bier construction on posets.
\end{abstract}

\date{September 2007}


\address{Fachbereich Mathematik, Universität Bremen, 28359 Bremen, Germany}

\email{jlehmann@math.uni-bremen.de}

\maketitle

\section{Introduction}

\emph{Nested set complexes} for semilattices were introduced by Feichtner
and Kozlov in their paper \cite{FK04} as a unifying framework for
the study of De Concini-Procesi models of subspace arrangements and
the resolution of singularities in toric varieties. Feichtner and
Müller considered the topology of those complexes (\cite{FM}), in
particular they prove that the nested set complex of any building set of a semilattice is
homotopy equivalent to the order complex of the semilattice without
its minimal element. These results found applications in the study
of complexes of trees (\cite{Fei05}) and $k$-trees (\cite{D}).

The \emph{Bier construction} was originally introduced by Thomas Bier
in 1992 (\cite{Bie92}) as a construction on abstract simplicial complexes;
more precisely, given an abstract simplicial complex $\mathcal{A}$,
the deleted join of $\mathcal{A}$ with the combinatorial Alexander
dual of $\mathcal{A}$ is another complex, the Bier sphere of $\mathcal{A}$.
A short proof that this construction actually results in a sphere
was given by De Longueville (\cite{Lon}). In 2004, Björner, Paffenholz,
Sjöstrand and Ziegler (\cite{BPSZ05}) reinterpreted the construction
in order-theoretic terms, by viewing an abstract simplicial complex
as an ideal in a Boolean lattice. Then the corresponding Bier poset
can be obtained as a subposet of the interval poset of the Boolean
lattice; that is, the poset consisting of only those intervals that cross the
ideal. This lends itself to immediate generalization, by considering
arbitrary bounded posets instead of a Boolean lattice. It turned out
that even in this general case the order complex of the Bier poset
is a subdivision of the order complex of the original poset. The complexes
that occur as intermediate steps of the subdivision are in general
not order complexes, as remarked in \cite{BPSZ05}.

Another view on the subject was taken by \v{C}uki\'{c} and Delucchi,
who in \cite{CD07} employed the theory of nested set complexes as a
framework for the study of the Bier construction. They found a new
proof for the result of Björner et al. for semilattices, exhibiting
the intermediate complexes as nested set complexes.

In this paper, we generalize the concept of nested sets to 
posets and exhibit the topological relationship between the arising
nested set complexes and the order complex of the underlying poset. A 
sufficient condition is given so that this relationship is actually a 
subdivision. Using these results, we generalize the proof method of 
\v{C}uki\'{c} and Delucchi to posets and obtain a new proof for the 
result of Björner et al. in full generality, showing that the intermediate 
complexes are still nested set complexes.

This paper is organized as follows. The terminology used is given
in Section~2. In Section~3, we extend the notion of nested set complexes
to posets, and prove our main Theorem~\ref{lem:expand-to-sd-or-cone}
about the topological behaviour of nested set complexes under extension
of the building set. In particular, under certain circumstances, subdivisions
take place, as it is generally the case when considering semilattices.
This allows to apply the framework to the treatment of the Bier construction
in Section~4, in the same way as in \cite{CD07}, to obtain a new
proof for the result of Björner et al. in full generality.

\begin{acknowledgement*}
I would like to thank Dmitry Kozlov for introducing me to the problem
and also him and Eva Maria Feichtner for the helpful discussions.
\end{acknowledgement*}

\section{Terminology}

In this paper, $P$ will denote a poset of finite length ($P$ will
be used for the underlying set interchangeably). For a general reference
on posets, see e.g. \cite{DP}. We will here use the following terminology:
Let $a,b\in P$ be elements of $P$, let $X$ be a subset of $P$. If $a\leq x$ 
for all $x\in X$, we write $a\leq X$; analogously $X\leq a$ means $x\leq a$ 
for all $x\in X$.

$P_{\leq X}$ denotes the set $\{a\in P:a\leq X\}$;
analogously we write $P_{\geq X}$ for $\{a\in P:a\geq X\}$. The set of \emph{upper bounds} of $X$ 
is $\ub X:=\min P_{\geq X}$; the set of \emph{lower bounds} of $X$ is
$\lb X:=\max P_{\leq X}$. If $P_{\geq X}$
has a least element $y$, so that $\ub X=\{y\}$, then $y$ is called the \emph{join} of $X$,
denoted by $\bigvee_{x\in X}x$ or simply by $\bigvee X$. Conversely, if $\lb X=\{y\}$ then
$y$ is called the \emph{meet} of $X$, denoted by $\bigwedge X$. For $\bigvee\{a,b\}$ and $\bigwedge\{a,b\}$,
the notations $a\vee b$ and $a\wedge b$ will be used, respectively.

We recall the difference between a poset and a (meet-)semilattice: In a semilattice
$P$, for any finite subset $X$ of $P$, either the join of $X$
exists or $\ub X$ is empty. But in a poset, sets of the form $\ub X$ with
$X\subset P,|\ub X|\geq 2$ can occur; these will be termed \emph{big cuts}.

If $P$ has a least element, this will be denoted by $\hat{0}$;
a greatest element will be denoted by $\hat{1}$. A \emph{bounded poset} is
a poset possessing $\hat{0}$ and $\hat{1}$; $\overline{P}$ means
$P\backslash\{\hat{0},\hat{1}\}$. For elements $x\leq y$ of $P$,
the \emph{interval} $[x,y]$ is defined as the poset with elements $z\in P$ where
$x\leq z\leq y$ and the order induced by $P$. An \emph{ideal} (order ideal
or down-set) $I$ of $P$ is a subset of $P$ with the property that
$x\in I$ and $y\leq x$ imply $y\in I$. In particular, if $P$ has a least element
$\hat{0}$, then every ideal contains $\hat{0}$.

Now recall the definition of a building set of a poset, as given in
\cite{FK04}:

\begin{defn}
Let $P$ be a poset with $\hat{0}$, let $G$ be a subset of $P_{>\hat{0}}$.
Denote by $F_{G}(x)$ the factors of $x$ in $G$, that is $\max G_{\leq x}=\max\{ g\in G:g\leq x\}$.
Then $G$ is a building set of $P$ if for any $x\in P$ there is
an isomorphism of posets \begin{eqnarray*}
\psi_{x}:\Pi_{i=1}^{t}[\hat{0},x_{i}] & \rightarrow & [\hat{0},x]\end{eqnarray*}
satisfying $\psi((\hat{0},\ldots,\hat{0},x_{i},\hat{0},\ldots,\hat{0}))=x_{i}$
for all $i$, where $\{ x_{1},x_{2},\ldots,x_{t}\}=F_{G}(x)$.
\end{defn}
\begin{figure}
\includegraphics{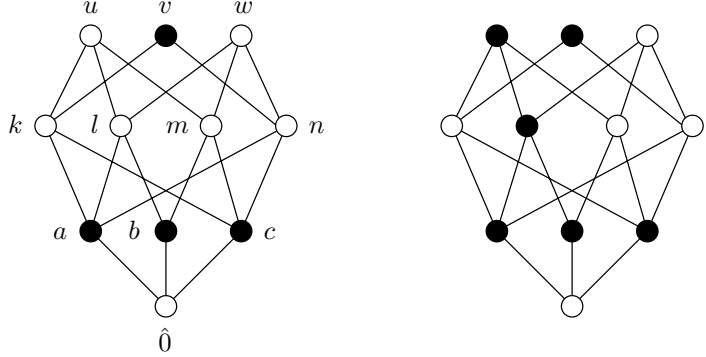}

\caption{\label{cap:buildingset}A poset with different building sets, marked
in black.}
\end{figure}
Figure~\ref{cap:buildingset} shows a poset with different building
sets.

For an abstract simplicial complex $\mathcal{A}$, the set of simplices
of $\mathcal{A}$ will be denoted with $\mathcal{A}$ as well. In
this paper, all abstract simplicial complexes are considered to contain
the empty set. If $X$ is the subset of the vertex set of an abstract
simplicial complex $\mathcal{A}$, then the \emph{subcomplex of $\mathcal{A}$
induced by $X$} has vertex set $X$ and set of simplices $\{\sigma\in\mathcal{A}:\sigma\subseteq X\}$.

In an abstract simplicial complex $\mathcal{A}$ with a face $\sigma$,
the \emph{stellar subdivision of $\mathcal{A}$ at $\sigma$} is an abstract
simplicial complex $\sd_{\mathcal{A}}\sigma$ with vertex set consisting
of the vertex set of $\mathcal{A}$ and an additional vertex $\hat{\sigma}$,
and with set of simplices defined as following: \begin{eqnarray*}
\sd_{\mathcal{A}}\sigma & := & \{\tau\in\mathcal{A}:\tau\not\supseteq\sigma\}\cup\{\tau\cup\{\hat{\sigma}\}:\tau\in\mathcal{A},\tau\not\supseteq\sigma,\tau\cup\sigma\in\mathcal{A}\}.\end{eqnarray*}
 For a subcomplex $\tilde{\mathcal{A}}$ of $\mathcal{A}$, the \emph{cone
over $\tilde{\mathcal{A}}$} is an abstract simplicial complex $\cone_{\mathcal{A}}\tilde{\mathcal{A}}$
with vertex set again consisting of the vertex set of $\mathcal{A}$
and an additional vertex $a$, and with set of simplices defined
as following: \begin{eqnarray*}
\cone_{\mathcal{A}}\tilde{\mathcal{A}} & := & \mathcal{A}\cup\{\tau\cup\{ a\}:\tau\in\tilde{\mathcal{A}}\}.\end{eqnarray*}

The \emph{order complex} $\Delta(P)$ of $P$ is an abstract simplicial complex
consisting of all chains of $P$.

\section{Nested set complexes and their behaviour under extension of the building
set}

The notion of a nested set introduced in \cite{FK04} can be extended to posets as it is.

\begin{figure}
\includegraphics{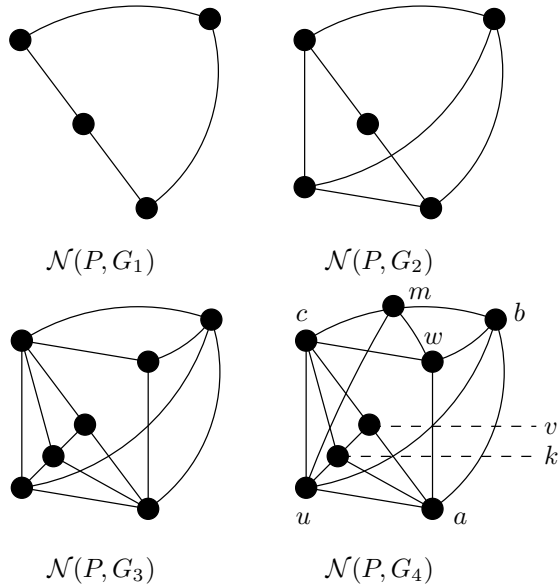}

\caption{\label{cap:nestedcomp}Nested set complexes of the poset $P$ shown
in Figure~\ref{cap:buildingset}, with the following building sets
(from left to right): $G_{1}=\{ a,b,c,v\}$, $G_{2}=\{ a,b,c,v,u\}$,
$G_{3}=\{ a,b,c,v,u,w,k\}$, $G_{4}=\{ a,b,c,v,u,w,k,m\}$. All triangles
are part of the corresponding complex.}
\end{figure}

\begin{defn}
Let $P$ be a poset of finite length with $\hat{0}$, let $G$ be
a building set of $P$. We call a finite subset $N\subset G$ nested
if for every incomparable subset $A\subset N$ with $|A|\geq2$ the
join of $A$ exists and $\bigvee A\notin G$.
\end{defn}
The nested sets in $G$ form an abstract simplicial complex, denoted
$\mathcal{N}(P,G)$ (see Figure~\ref{cap:nestedcomp} for some examples). 
Note that every singleton subset of $G$ is nested in $G$, thus the
vertices of $\mathcal{N}(P,G)$ correspond to the elements of $G$. 
Moreover, extending the building set has topological significance
for the nested set complex:

\begin{thm}
\label{lem:expand-to-sd-or-cone}Let $P$ be a poset of finite length
with $\hat{0}$, let $G$ be a building set of $P$. Let $x\in\max P\backslash G$.
Then $\tilde{G}:=G\cup\{ x\}$ is a building set of $P$ and 

(1) $\mathcal{N}(P,\tilde{G})\cong\sd_{\mathcal{N}(P,G)}F_{G}(x)$
if $F_{G}(x)=\{ x_{1},\ldots,x_{t}\}$ is a face of $\mathcal{N}(P,G)$,

(2) $\mathcal{N}(P,\tilde{G})\cong\cone_{\mathcal{N}(P,G)}\mathcal{C}$
otherwise, where $\mathcal{C}$ is the subcomplex of $\mathcal{N}(P,G)$
induced by $C_{G}(x):=G_{<x}\cup G_{>x}$.
\end{thm}

\begin{proof}
$\tilde{G}$ is clearly a building set of $P$.

Note that $F_{G}(x)$ is finite, since any interval in a poset of
finite length cannot be isomorphic to an infinite product of non-trivial
posets.

For (1), we need to see that \begin{eqnarray*}
\mathcal{N}(P,\tilde{G}) & = & \{ N\in\mathcal{N}(P,G):N\not\supset F_{G}(x)\}\cup\\
 &  & \{ N\cup\{ x\}:N\in\mathcal{N}(P,G),N\not\supset F_{G}(x),N\cup F_{G}(x)\in\mathcal{N}(P,G)\}\\
 & \cong & \sd_{\mathcal{N}(P,G)}F_{G}(x),\end{eqnarray*}
for (2), we need to see that\begin{eqnarray*}
\mathcal{N}(P,\tilde{G}) & = & \mathcal{N}(P,G)\cup\{ N\cup\{ x\}:N\in\mathcal{N}(P,G),N\subset C_{G}(x)\}\\
 & \cong & \cone_{\mathcal{N}(P,G)}\mathcal{C}.\end{eqnarray*}

To this end, we show three equivalences:\[
\tag{a}\mathcal{N}(P,G)\cap\mathcal{N}(P,\tilde{G})=\{ N\in\mathcal{N}(P,G):N\not\supset F_{G}(x)\}=:\tilde{\mathcal{N}}.\]
 If $F_{G}(x)\in\mathcal{N}(P,G)$, then \[
\tag{b}\mathcal{N}(P,\tilde{G})\backslash\mathcal{N}(P,G)=\{ N\cup\{ x\}:N\in\tilde{\mathcal{N}},N\cup F_{G}(x)\in\mathcal{N}(P,G)\},\]
 if $F_{G}(x)\notin\mathcal{N}(P,G)$, then \[
\tag{c}\mathcal{N}(P,\tilde{G})\backslash\mathcal{N}(P,G)=\{ N\cup\{ x\}:N\in\mathcal{N}(P,G),N\subset C_{G}(x)\}.\]
 In Case 1, (a) and (b) give the result above; in Case 2, (a) and
(c) are needed. Note that in Case 2 no nested set can contain $F_{G}(x)$,
so $\tilde{\mathcal{N}}=\mathcal{N}(P,G)$ in that case.

In this proof let $\psi:=\psi_{x}$.

(a) Let $N$ be nested in $G$ and not nested in $\tilde{G}$. Then
an incomparable subset $A=\{ a_{1},\ldots,a_{s}\}$ of $N$ exists
with $\bigvee A=x$, so $A\subset[\hat{0},x]$. Let $a_{i}=\psi(a_{i1},a_{i2},\ldots,a_{it})$,
then since for all $a_{i}$ exists $x_{j_{i}}$ such that $a_{i}\leq x_{j_{i}}$,
we have $a_{ij}=\hat{0}$ for all $j\neq j_{i}$. So consider $A_{j}=\{ a_{i}\in A:j_{i}=j\}$;
since this is an incomparable subset of $A$ which is nested in $G$,
the join of $A_{j}$ exists. Since the join of $A$ exists, it must coincide
with the join of $\psi^{-1}(A)$ in $\Pi_{i}[\hat{0},x_{i}]$, same
for each $A_{j}$. Hence $\bigvee A_{j}=\psi(\bigvee\psi^{-1}(A_{j}))=x_{j}$
which is in $G$. So $|A_{j}|\leq1$ holds; $A_{j}=\emptyset$ for any
$j$ implies $\bigvee A\neq x$, thus $A_{j}=\{ x_{j}\}$ for all $j$,
so $A=F_{G}(x)\subset N$.

\begin{sloppypar}Let $N$ be nested in $G$ and $F_{G}(x)\subset N$.
$F_{G}(x)$ is incomparable, and $\bigvee F_{G}(x)$ exists since
$N$ is nested in $G$, so $\bigvee F_{G}(x)=\psi(\bigvee_{x_{i}\in F_{G}(x)}\psi^{-1}(x_{i}))=\psi(x_{1},x_{2},\ldots,x_{t})=x$,
which is not in $\tilde{G}$. Thus $N$ is not nested in $\tilde{G}$.\end{sloppypar}

(b) Let $F_{G}(x)$ be nested in $G$, $N\subset\tilde{G}$ containing
$x$, $N\backslash\{ x\}\not\supset F_{G}(x)$ and $(N\backslash\{ x\})\cup F_{G}(x)$
be nested in $G$. Then $N$ is not nested in $G$ since $N\not\subset G$.
Since $N\backslash\{ x\}$ is nested in $G$ and does not contain
$F_{G}(x)$, $N\backslash\{ x\}$ is nested in $\tilde{G}$ by (a).
So only sets $A\subset N$ incomparable with $x\in A$ and $|A|\geq2$
have to be investigated further. Note that for all $a\in A$ there
is no $x_{i}$ with $a\leq x_{i}$, since otherwise $a\leq x$ follows
in contradiction to $A$ being incomparable. Let $\tilde{A}=(A\backslash\{ x\})\cup\{ x_{i}\in F_{G}(x):x_{i}$ is incomparable to all $a\in A\backslash\{ x\}\}$.
Assume $|\tilde{A}|=1$, that means that $A=\{ a,x\}$ and $a>x_{i}$
for all $x_{i}$, hence $a\geq\bigvee F_{G}(x)=x$, so $A$ is not
incomparable. Hence $|\tilde{A}|\geq2$, and since $\tilde{A}\subset(N\backslash\{ x\})\cup F_{G}(x)$
and $\tilde{A}$ is incomparable, $\bigvee\tilde{A}$ exists and is
not in $G$. But $\bigvee\tilde{A}\geq F_{G}(x)$, so $\bigvee\tilde{A}\geq x$
with equality only if for all $a\in A$ there is $x_{i}$ with $x\geq a>x_{i}$
in contradiction to $\tilde{A}\subset G$. Hence $\bigvee\tilde{A}>x$,
and therefore by the choice of $x$ the join of $\tilde{A}$ is contained in $G$.

So there exists no $A\subset N$ incomparable with $x\in A$ and $|A|\geq2$,
so $N$ is nested in $\tilde{G}$.

Now let $F_{G}(x)$ be nested in $G$ and $N$ be nested in $\tilde{G}$,
but not nested in $G$. As the existence of the join is independent
of the considered building set, all nested sets in $\tilde{G}$ not
containing $x$ are nested in $G$. So $x\in N$ and $N\backslash\{ x\}$
is nested in $G$ and in $\tilde{G}$. Thus by (a), $N\backslash\{ x\}\not\supset F_{G}(x)$.

Consider $B\subset(N\backslash\{ x\})\cup F_{G}(x)$ incomparable,
$|B|\geq2$, containing $x_{s}\in F_{G}(x)\backslash N$. Let $B_{N}=B\backslash F_{G}(x)\subset N\backslash\{ x\}$,
$B_{F}=B\cap F_{G}(x)\subset F_{G}(x)$. Since $F_{G}(x)$ is nested
in $G$, $\bigvee B_{F}$ exists and is in $[\hat{0},x]\backslash G$.
If $B_{N}$ is empty, then $B=B_{F}$, and $\bigvee B\notin G$ exists. 

So assume $|B_{N}|\geq1$ in the following. As $N\supset B_{N}\cup\{ x\}$
is nested in $\tilde{G}$, either $\bigvee B_{N}\cup\{ x\}\notin\tilde{G}$
exists, which is not possible since $\bigvee B_{N}\cup\{ x\}\geq x$
which by the choice of $x$ is in $\tilde{G}$, or $B_{N}\cup\{ x\}$
is not incomparable. Assuming the existence of $b\in B_{N}$ with
$b>x$ implies $b>x>x_{s}$, so $B$ would not have been incomparable.
Let $\tilde{B}=\{ x\}\cup\{ b\in B_{N}:b\textnormal{ is incomparable to }x\}$;
$\tilde{B}\subset N$ is incomparable. Assuming $|\tilde{B}|\geq2$
we obtain a contradiction to $\bigvee\tilde{B}\notin\tilde{G}$ as
above and thus to $N$ nested in $\tilde{G}$. Thus for all $b\in B_{N}$
exists $x_{i}$ with $b<x_{i}$, so $B_{N}\subset[\hat{0},x]$. Since
$B$ is incomparable, these $x_{i}$ are all not in $B_{F}$, so there
is w.l.o.g. a partition $\{[f],[t]\backslash[f]\}$ of $[t]$, where
$[f]=\{ i:x_{i}\in B_{F}\}$, and $\psi^{-1}(B_{N})\subset\Pi_{i\in[f]}\{\hat{0}\}\times\Pi_{i\in[t]\backslash[f]}[\hat{0},x_{i}]$.
So $\psi^{-1}(\bigvee B_{N})=(\hat{0},\ldots,\hat{0},y_{f+1},y_{f+2},\dots,y_{t})$
and $\psi^{-1}(\bigvee B_{F})=(y_{1},y_{2},\ldots,y_{f},\hat{0},\ldots,\hat{0})$
(both joins exist because $F_{G}(x),N\ \backslash\{ x\}$ are nested
in $G$), so $\bigvee B=\bigvee B_{F}\vee\bigvee B_{N}=\psi(y_{1},y_{2},\ldots,y_{t})$
exists. Since $B_{F},B_{N}\neq\emptyset$, at least two $y_{i}\neq\hat{0}$,
so by definition of $F_{G}(x)$ the join of $B$ is not in $G$.
Thus $(N\backslash\{ x\})\cup F_{G}(x)$ is nested in $G$.

(c) Let $F_{G}(x)$ be not nested in $G$ and $N\in\mathcal{N}(P,\tilde{G})\backslash\mathcal{N}(P,G)$.
As in (b), $x\in N$ and $N\backslash\{ x\}$ nested in $G$. Let
$n\in N\backslash\{ x\}$. If $\{ n,x\}\subset N$ is incomparable,
then $n\vee x$ exists and is not in $\tilde{G}$, which is impossible
by the choice of $x$ since $n\vee x\geq x$. Thus $n\in C_{G}(x)$.

Let $F_{G}(x)$ be not nested in $G$, let $N\subset\tilde{G}$ with
$x\in N$ and let $N\backslash\{ x\}$ be a subset of $C_{G}(x)$
that is nested in $G$. Since $x\in N$, $N$ is not a subset of $G$,
so not nested in $G$. Let $A\subset N$ be incomparable with $|A|\geq2$,
so $A\subset N\backslash\{ x\}$, so $\bigvee A\notin G$ exists.
Hence by choice of $x$, $A\subset G_{<x}$, that is for all $a\in A$
exists $x_{i}\in F_{G}(x)$ with $a\leq x_{i}$. As $x\in\ub F_{G}(x)$
(by isomorphism of $[\hat{0},x]$ to the product of the intervals
$[\hat{0},x_{i}]$) and $F_{G}(x)$ is incomparable but not nested
in $G$, $\bigvee F_{G}(x)$ does not exist, so there exists $y\in\ub F_{G}(x)$,
$y\neq x$. Since $y\geq A$, assuming $\bigvee A=x$ implies $y\geq x$,
a contradiction to the choice of $y$. Hence $\bigvee A$ is even not
in $\tilde{G}=G\cup\{ x\}$, so $N$ is nested in $\tilde{G}$.
\end{proof}

\begin{figure}
\includegraphics{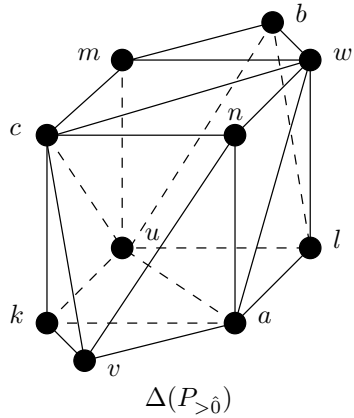}

\caption{\label{cap:ordercomp}The reduced order complex of the poset $P$ shown in
Figure~\ref{cap:buildingset}, identical to the nested set complex
for the building set $P_{>\hat{0}}$. Again, all triangles
are part of the complex.}
\end{figure}

Note that for the maximal building set $P_{>\hat{0}}$ the nested
set complex $\mathcal{N}(P,P_{>\hat{0}})$ coincides with the order
complex $\Delta(P_{>\hat{0}})$ of $P_{>\hat{0}}$ by definition.
Thus, by successively expanding a building set of $P$ as in the preceding
theorem, we get a sequence of nested set complexes, eventually arriving
at the order complex of the poset. Some steps of this process are shown in
Figures~\ref{cap:nestedcomp} and \ref{cap:ordercomp}. Now it turns
out that by considering the big cuts of a poset, a guarantee can be
given for some building sets to only yield stellar subdivisions in
every expansion step.

\begin{lem}
\label{lem:all-cuts->factors-nested}If $G$ is a building set of
$P$ with $\ub A\subseteq G$ for all $A\subseteq P$ with $|\ub A|\geq2$,
then $F_{G}(x)\in\mathcal{N}(P,G)$ for all $x\in P$.
\end{lem}
\begin{proof}
Let $x\in P$; if $x\in G$ then $F_{G}(x)=\{ x\}$ is nested in $G$.
If $x\notin G$, then $F_{G}(x)=\{ x_{1},x_{2},\ldots,x_{t}\}$, let $\psi:=\psi_{x}$. Let
$A\subseteq F_{G}(x)$ with $|A|\geq2$. Then $z=\psi(z_{1},\ldots,z_{t})\in\ub A$,
where $z_{i}=x_{i}$ if $x_{i}\in A$ and $z_{i}=\hat{0}$ otherwise.
Since we have $x_{i}<z\leq x$ for all $x_{i}\in A$, $z\notin G$
holds by the definition of $F_{G}(x)$. So $\ub A=\{ z\}$, or in
other words, the join of $A$ exists (and equals $z$). Thus $F_{G}(x)$
is nested in $G$.
\end{proof}
\begin{thm}
For any building set $G$ of $P$ fulfilling the condition in Lemma~\ref{lem:all-cuts->factors-nested},
$\Delta(P_{>\hat{0}})$ is a subdivision of $\mathcal{N}(P,G)$.
\end{thm}
\begin{proof}
Stepwise expanding $G$ to the maximal building set $P_{>\hat{0}}$
along a linear extension of $P_{>\hat{0}}\backslash G$ yields a stellar
subdivision of the corresponding nested set complex in each step by
Lemma~\ref{lem:all-cuts->factors-nested} and Theorem~\ref{lem:expand-to-sd-or-cone},
so $\mathcal{N}(P,G)$ subdivides to $\mathcal{N}(P,P_{>\hat{0}})=\Delta(P_{>\hat{0}})$.
\end{proof}

\section{Application to Bier posets}

The condition given in Lemma~\ref{lem:all-cuts->factors-nested}
is not at all necessary as we will see now. The special structure
of Bier posets allows to determine easily that for a certain building
set all extensions induce only stellar subdivisions of the corresponding
nested set complexes, despite the occurence of possibly many big cuts
outside of this building set.

But first of all we will recall the definition of a Bier poset.

\begin{defn}
Given a poset $P$ of finite length with $\hat{0}$ and a proper ideal
$I$ of $P$, the Bier poset $\Bier(P,I)$ is a poset with

\textbullet\, elements $\{\hat{1}\}\cup\{[x,y]\textnormal{ intervals of }P:x\in I,y\notin I\}$,

\textbullet\, order relation $[x,y]\leq[v,w]$ iff $x\leq v<w\leq y$
and $[x,y]\leq\hat{1}$ for all $I\ni x<y\notin I$.
\end{defn}
We now present a new structural proof using nested sets of the following 
result of Björner, Paffenholz, Sjöstrand and Ziegler in its full generality:

\begin{thm}
\label{thm:bpsz-bier-is-edge-sd}\cite[Thm. 2.2]{BPSZ05} Let $P$
be a bounded poset of finite length with a proper ideal $I$. Then the order
complex of $\overline{\Bier(P,I)}$ is obtained from the order complex
of $\bar{P}$ by stellar subdivision of the edges $\{ x,y\}$ with
$x<y$, $x\in I_{>\hat{0}}$, $y\in\bar{P}\backslash I$ in order
of increasing length of the corresponding intervals.
\end{thm}
The argumentation follows \cite{CD07}, starting with a building set
for $\Bier(P,I)_{<\hat{1}}$:

\begin{prop}
\cite[Prop. 2.1]{CD07} For any bounded poset $P$ of finite length with proper
ideal $I$, \begin{eqnarray*}
G & = & \{[x,\hat{1}],[\hat{0},y]:x\in I_{>\hat{0}},y\in\bar{P}\backslash I\}\end{eqnarray*}
is a building set of $\Bier(P,I)_{<\hat{1}}$.
\end{prop}
The proof of the above proposition in \cite{CD07} does not use the
lattice property or the finiteness and thus remains valid in the poset case. This building
set of the Bier poset is very well behaved: Let $x\in I_{>\hat{0}},y\in\bar{P}\backslash I$,
let $[a,b]\in\overline{\Bier(P,I)}_{\geq\{[x,\hat{1}],[\hat{0},y]\}}$.
Then we have $x\leq a<b\leq\hat{1}$ and $\hat{0}\leq a<b\leq y$,
so $[x,y]\leq[a,b]$. Hence $[x,\hat{1}]\vee[\hat{0},y]=[x,y]$.

This enables us to follow \v{C}uki\'{c} and Delucchi further; we
find that their characterization of the nested sets in $G$ remains
valid as well:

\begin{lem}
\label{lem:nested-sets-in-G-Bier}\cite[Lem. 2.2]{CD07} With $G$
as above, a set $N\subset G$ is nested in $G$ iff
\begin{enumerate}
\item If $[\hat{0},y_{1}]$, $[\hat{0},y_{2}]\in N$, then $y_{1}$
and $y_{2}$ are comparable.
\item If $[x_{1},\hat{1}],[x_{2},\hat{1}]\in N$, then $x_{1}$ and
$x_{2}$ are comparable.
\item If $[x,\hat{1}],[\hat{0},y]\in N$, then $x<y$.
\end{enumerate}
\end{lem}
\begin{proof}
Let $N$ be nested in $G$. Assume that (i) is not true, that is, $[\hat{0},y_{1}],[\hat{0},y_{2}]\in N$
exist with $y_{1},y_{2}$ incomparable. Then $[\hat{0},y_{1}],[\hat{0},y_{2}]$
are incomparable as well, so since $N$ is nested, $[\hat{0},y_{1}]\vee[\hat{0},y_{2}]$
exists and is not an element of $G$.

Let $S=\lb\{ y_{1},y_{2}\}$. If $p\in S\cap I$ exists, then for
all $[a,b]\in\overline{\Bier(P,I)}_{\geq\{[\hat{0},y_{1}],[\hat{0},y_{2}]\}}$,
we have $\hat{0}\leq a<b\leq y_{1},y_{2}$, thus $b\leq p\in I$ holds,
so $b\in I$, a contradiction to $[a,b]\in\Bier(P,I)$. Hence $[\hat{0},y_{1}]\vee[\hat{0},y_{2}]$
does not exist. So in the following let $S\cap I$ be empty.

\begin{sloppypar}Case 1: If $y_{1}\wedge y_{2}$ exists (and, as
noted above, then is in $P\backslash I$), consider $[a,b]\in\overline{\Bier(P,I)}_{\geq\{[\hat{0},y_{1}],[\hat{0},y_{2}]\}}$,
which means $\hat{0}\leq a<b\leq y_{1}$ and $\hat{0}\leq a<b\leq y_{2}$,
and thus $[\hat{0},y_{1}\wedge y_{2}]\leq[a,b]$. Hence $[\hat{0},y_{1}]\vee[\hat{0},y_{2}]=[\hat{0},y_{1}\wedge y_{2}]\in G$,
a contradiction.\end{sloppypar}

Case 2: If $y_{1}\wedge y_{2}$ does not exist, this means that $|S|\geq2$.
Observe that $[\hat{0},s]\geq[\hat{0},y_{i}]$ for $i=1,2$ and all $s\in S$. Let $[p,q]\geq[\hat{0},y_{i}]$
for $i=1,2$; that is, $\hat{0}\leq p<q\leq y_{1},y_{2}$. By definition
of $S$, there exists $s\in S$ with $\hat{0}\leq p<q\leq s$, so
we have $[p,q]\geq[\hat{0},s]$. Thus, $[\hat{0},s]\in\hat{S}:=\ub\{[\hat{0},y_{1}],[\hat{0},y_{2}]\}$ for
all $s\in S$, so $|\hat{S}|\geq2$ meaning that $[\hat{0},y_{1}]\vee[\hat{0},y_{2}]$
does not exist.

Analogously we obtain (ii).

For (iii), we see that $[x,\hat{1}],[\hat{0},y]$ are incomparable
(since $x\neq\hat{0},y\neq\hat{1}$), so since $N$ is nested, $[x,\hat{1}]\vee[\hat{0},y]=[p,q]\in\Bier(P,I)_{<\hat{1}}$.
So $x\leq p<q\leq\hat{1}$ and $\hat{0}\leq p<q\leq y$ hold, implying
$x<y$.

Conversely, consider a set $N\subset G$, fulfilling conditions (i)-(iii).
Let \[A=\{[x_{1},\hat{1}],\ldots,[x_{s},\hat{1}],[\hat{0},y_{1}],\ldots,[\hat{0},y_{t}]\}\]
 be an incomparable subset of $N$, with $|A|\geq2$. By conditions
(i) and (ii), $s=t=1$, and as explained above, $[x_{1},\hat{1}]\vee[\hat{0},y_{1}]=[x_{1},y_{1}]\notin G$
since $x_{1}\neq\hat{0},y_{1}\neq\hat{1}$.
\end{proof}
So the nested sets in this particular building set $G$ of a Bier
poset coincide with those in the lattice case, which allows to follow
\v{C}uki\'{c} and Delucchi further, thus obtaining

\begin{prop}
\cite[Prop. 2.3]{CD07} For $P,I$ and $G$ as above, $\mathcal{N}(\Bier(P,I)_{<\hat{1}},G)=\Delta(\bar{P})$.
\end{prop}
\begin{proof}
As in \cite{CD07}, by Lemma~\ref{lem:nested-sets-in-G-Bier} all
nested sets $N$ in $G$ are of the form $N=\{[x_{1},\hat{1}],\ldots,[x_{s},\hat{1}],[\hat{0},y_{1}],\ldots,[\hat{0},y_{t}]\}$,
where $x_{1}<x_{2}<\ldots<x_{s}<y_{1}<y_{2}<\ldots<y_{t}$, and $f:\mathcal{N}(\Bier(P,I)_{<\hat{1}},G)\rightarrow\Delta(\bar{P})$,
mapping a nested set to its underlying chain, is an order-preserving
bijection.
\end{proof}
Now the good behavior of $G$ noted above comes into play once again. Considering
any element $[x,y]\in\overline{\Bier(P,I)}\backslash G$, the factors
of $[x,y]$ in $G$ are $F_{G}([x,y])=\{[x,\hat{1}],[\hat{0},y]\}$,
which is an incomparable set with join $[x,y]$. Thus $F_{G}([x,y])$
is a nested set in $G$, but also in all building sets $\tilde{G}$
resulting from repeated application of Theorem~\ref{lem:expand-to-sd-or-cone},
starting with the building set $G$. Hence in all applications of
Theorem~\ref{lem:expand-to-sd-or-cone} only stellar subdivisions
occur; more precisely we see that $\Delta(\overline{\Bier(P,I)}$
is obtained from $\mathcal{N}(\Bier(P,I)_{<\hat{1}},G)=\Delta(\bar{P})$
by stellar subdivisions of all edges $f(\{[x,\hat{1}],[\hat{0},y]\})=\{ x<y\}$,
in order of increasing length $\ell([x,y])$ in $P$ (since an interval of length
1 in $P$ is maximal in $\Bier(P,I)_{<\hat{1}}$, covering there the intervals of length 2
and so on). This concludes the
proof of Theorem~\ref{thm:bpsz-bier-is-edge-sd} via nested set complexes.

\end{document}